\documentclass[a4paper,11pt]{article}
\usepackage{amsmath,amsthm,amssymb}

\usepackage[utf8]{inputenc}
\usepackage[english]{babel}
\usepackage{hyperref}
\hypersetup{
    colorlinks = true,
    linkcolor = {black},
    citecolor = {black}
}

\usepackage{enumitem}
\usepackage{tikz}
\usetikzlibrary{calc,decorations.pathmorphing,decorations.text,fixedpointarithmetic}
\usepackage{subfigure}
\usepackage{refcount}
\usepackage[hmargin=3cm,vmargin=3cm]{geometry}
\usepackage{graphicx,color,xcolor}    % for latex
\usepackage{epsfig}
\usepackage{caption}
\usepackage{tikz}
\usetikzlibrary{calc}
\usetikzlibrary{backgrounds, shapes.geometric}
\usepackage{rotating}
\usetikzlibrary{positioning}
\usepackage{pgffor}
\usepackage{authblk}

\newtheorem{theorem}{Theorem}[section]

\newtheorem{conjecture}[theorem]{Conjecture}

\newcommand{\NN}{\mathbb{N}}

\def\COMMENT#1{}
\let\COMMENT=\footnote% COMMENT OUT for clean output

\usepackage{caption}

\begin{document}

\title{A bounded diameter strengthening of K\H onig's Theorem}
\author[1]{Louis DeBiasio\thanks{Research supported in part by NSF grant DMS-1954170.}}
\author[2]{Ant\'onio Gir{\~a}o\thanks{Research supported by ERC Advanced Grant no. 883810.}}
\author[3]{Penny Haxell\thanks{Partially supported by NSERC}} 
\author[4]{Maya Stein\thanks{Supported by FONDECYT Regular Grant 1221905 and by ANID Basal Grant CMM FB210005.
}}

\affil[1]{\small Department of Mathematics, Miami University, Oxford, OH. 
E-mail: debiasld@miamioh.edu}
\affil[2]{{Mathematical Institute, University of Oxford, Andrew Wiles Building, Radcliffe Observatory Quarter, Woodstock Road, Oxford, UK.}
E-mail: girao@maths.ox.ac.uk}
\affil[3]{\small Department of Combinatorics and Optimization,
          University of Waterloo, Waterloo ON Canada N2L 3G1. 
          E-mail: pehaxell@uwaterloo.ca
          }
\affil[4]{Department of Mathematical Engineering and Center for Mathematical Modeling (CNRS IRL2807), University of Chile. 
E-mail: mstein@dim.uchile.cl
}

\date{\today}
%\date{\bf {Jan 31, 2001}}

\maketitle

\begin{abstract}
K\H onig's theorem says that the vertex cover number of every bipartite graph is at most its matching number (in fact they are equal since, trivially, the matching number is at most the vertex cover number). An equivalent formulation of K\H onig's theorem is that in every $2$-colouring of the edges of a graph $G$, the number of monochromatic components needed to cover the vertex set of $G$ is at most the independence number of $G$.

We prove the following strengthening of K\H onig's theorem:  In every $2$-colouring of the edges of a graph $G$, the number of monochromatic subgraphs of bounded diameter needed to cover the vertex set of $G$ is at most the independence number of $G$.
\end{abstract}

\section{Introduction}

Let $G$ be a (hyper)graph.
A \emph{matching} in $G$ is a set of pairwise disjoint edges.  A \emph{vertex cover} of the edges of $G$ is a set of vertices $S$ such that every edge of $G$ is incident with a vertex from~$S$.  We denote the size of a largest matching in $G$ by $\nu(G)$ and the size of a minimum vertex cover of $G$ by $\tau(G)$.  Note that for every (hyper)graph $G$ we have $\nu(G)\leq \tau(G)$, since a minimum vertex cover must contain at least one vertex from each edge in a maximum matching.  

The following theorem of K\H onig from 1931 is one of the foundational results in graph theory. 

\begin{theorem}[K\H onig \cite{Kon}]\label{thm:K1}
For every bipartite graph $G$, we have $\tau(G)\leq \nu(G)$.
\end{theorem}

While it has been observed many times before (see \cite{Gy}), it is perhaps less well-known that Theorem \ref{thm:K1} is equivalent\footnote{The equivalence can be seen by viewing a bipartite graph as the intersection graph between red and blue components (see \cite{GySurv1} or \cite{DKMS} for a more detailed explanation).} to the following statement about arbitrary 2-colourings of the edges of a graph, in which each edge receives at least one of the colours. It relates the smallest number of monochromatic components in such a colouring that are needed to cover $V(G)$ to the \emph{independence number} $\alpha(G)$, i.e. the size of a largest independent set of vertices in $G$. Here the term \emph{cover} is used to mean that the union of the vertex sets of the components is the whole set $V(G)$.

\begin{theorem}[K\H onig \cite{Kon}]\label{thm:K2}
For every graph $G$ and every 2-colouring of the edges of $G$, there exists a set of at most $\alpha(G)$ monochromatic components that cover $V(G)$.
\end{theorem}

Note that when $G$ is a complete graph, Theorem~\ref{thm:K2} is equivalent to saying that a graph or its complement is connected, a fact sometimes referred to as a ``first exercise in graph theory" or ``a remark of Erd\H{o}s and Rado''.  Another version of this fact, normally attributed to ``folklore" (see \cite[Theorem 2.1.11]{W}) is a strengthening in terms of diameter. The \emph{diameter} of a graph $H$ is the smallest $t$ such that every pair of vertices in $H$ are joined by a path of length at most $t$.

\begin{theorem}[Folklore]\label{thm:folk}
In every red-blue colouring of the edges of a complete graph, either the red graph has diameter at most 2, the blue graph has diameter at most 2, or both graphs have diameter exactly 3.  
\end{theorem}

This raises the question of whether one can generalize Theorem \ref{thm:folk} to obtain an analogous strengthening of K\H onig's theorem.  This conjecture was indeed made in \cite{DKMS}. We will say that a set $S$ of monochromatic subgraphs in an edge-coloured graph $G$ is a \emph{covering set} for $G$ if $V(G)$ is contained in the union of the vertex sets of the members of $S$. 
\begin{conjecture}[DeBiasio, Kamel, McCourt, Sheats \cite{DKMS}]\label{con:DKMS}
There exists a function $f:\NN\to \NN$ such that for every graph $G$ and every 2-colouring of the edges of $G$, there exists a covering set for $G$ of at most $\alpha(G)$ monochromatic subgraphs, each of which has diameter at most $f(\alpha(G))$.
\end{conjecture}
Observe that the monochromatic subgraphs in this conjecture might be contained in monochromatic components of larger diameter.
Theorem \ref{thm:folk} shows that $f(1)=3$, and in \cite{DKMS} the authors prove that $f(2)\leq 6$.  
Here we verify Conjecture~\ref{con:DKMS} by proving the following.

\begin{theorem}\label{thm:poly}
For every graph $G$ and every 2-colouring of the edges of $G$, there
exists a covering set for $G$ of at most $\alpha:=\alpha(G)$ monochromatic
subgraphs, each of which has diameter at most $8\alpha^2+12\alpha+6$.
\end{theorem}

We remark that for simplicity we do not attempt to get the best possible quadratic upper bound on the diameter (as we have no reason to believe that the optimal upper bound should be quadratic). 
In terms of lower bounds, currently we do not know any bound for the function $f$ besides the trivial $f(\alpha)\geq 3$ from Theorem~\ref{thm:folk}.

\section{Proof of Theorem \ref{thm:poly}}\label{sec:polyproof}
Let $G$ be a graph with a fixed red-blue colouring of its edges.

First we establish some notation.
For a vertex $v$ and integer $d\geq 0$, we denote by $N^R_d(v)$ the set of
vertices at red-distance $d$ from $v$, and by $N_{\leq d}^R(v)$ the set $\bigcup_{i\leq d}N^R_i(v)$. We write $N^R(v)$ for the set of neighbours of $v$ via red edges. We extend these definitions to subsets $S\subseteq V(G)$ by setting $N_{\leq d}^R(S)=\bigcup_{v\in S}N_{\leq d}^R(v)$ and $N^R(S)=\bigcup_{v\in S}N^R(v)$. We make the analogous definitions for the blue edges. Thus in particular $N(S)=N^R(S)\cup N^B(S)$.

We prove Theorem~\ref{thm:poly} by induction on $\alpha:=\alpha(G)$. If $\alpha=1$ the result is
true by Theorem~\ref{thm:folk}, so assume $\alpha\geq 2$ and that the theorem holds for all smaller values. Set $f(\alpha)=8\alpha^2+12\alpha+6$. 
If every monochromatic component in $G$ has diameter less than $f(\alpha)$
then we are done by Theorem~\ref{thm:K2}, so let us assume
without loss 
of generality that there exists a red component of diameter at least $f(\alpha)$. Then there exists a vertex $z$ and another vertex of red-distance from $z$ that is finite but at least $f(\alpha)$, so 
in particular the set $N^R_{d}(z)$ is  nonempty for all
$0\leq d\leq f(\alpha)$. For simplicity of notation we set 
$$r:=f(\alpha)/2-1=4\alpha^2+6\alpha+2$$ and $U:=N^R_{\leq r}(z).$

Suppose  $U$  contains  distinct nonadjacent vertices $x$ and $y$ that
are joined by a blue path of length less than $r$. 
Then $\{x,y\}\cup N(x)\cup N(y)\subseteq N_{\leq r+1}^R(z)\cup N_{\leq
  r}^B(x)$. Hence by 
removing $N_{\leq r+1}^R(z)\cup N_{\leq  r}^B(x)$ from $G$ we obtain a
graph $G'$ such that  
$\alpha(G')\leq\alpha-2$. By the induction hypothesis $G'$ has a
covering set of at most $\alpha-2$ monochromatic subgraphs, each of
diameter at most $f(\alpha-2)<f(\alpha)$. Together with
$N_{\leq r+1}^R(z)$ and $N_{\leq  r}^B(x)$, each of which has diameter at most $2(r+1)=f(\alpha)$, they form the required covering set for
$G$. So from now on we assume no such pair $\{x,y\}$ exists in~$U$; that is,
\begin{equation}\label{eq:blue}
\text{the blue-distance between any two distinct nonadjacent $x,y\in U$ is at least $r$.}
\end{equation}

Let $I$ be a maximal independent set of $G[U]$. 
Then each $u\in U\setminus(N^R_{\leq 1}(I))$ is
blue-adjacent to some vertex $a\in I$. By \eqref{eq:blue}, this
vertex is unique, and we call $a$ the
\emph{label} $\ell(u)$ of $u$.  Note
that 
\begin{equation}\label{eq:differentlabels}
\text{two vertices with different labels cannot be blue-adjacent,}
\end{equation}
as otherwise a blue path of length 3 would exist between their
labels, contradicting~\eqref{eq:blue}. Moreover, \begin{equation}\label{eq:samelabels}
\text{every pair of vertices with the same label must be adjacent
(in some colour),}
\end{equation} or else we again contradict \eqref{eq:blue}.

For a non-negative integer $t\leq\alpha$, we say that
a set $S$ of $t$ vertices is $t$-\emph{good} 
%for $I$
if the following hold.
\begin{enumerate}[label = {\bfseries{(G\arabic{enumi}})}]
\item\label{g1} $S\subseteq N^R_{\leq r-\alpha-2+t}(z)$,
\item\label{g2} every $s\in S$ is at red-distance at least $2\alpha-2t+3$ from
    every vertex in $I\cup (S\setminus\{s\})$,
    %\item the set $N_I^B(S)$ has size $t$, 
\item\label{g3} $\ell(s)\not=\ell(s')$ for each $s\not= s'\in S$.
\end{enumerate}
 Observe that by \ref{g1} and \ref{g2} all 
vertices in $N^R_{\leq 1}(S)$ have labels, so the labels in \ref{g3} are well-defined. Also, \ref{g2} implies
$S\cap I=\emptyset$. By \ref{g2}, \ref{g3} and \eqref{eq:differentlabels} we know that
$S$ is independent.

Let $t$ be largest such that a $t$-good set $S$ exists.  
We claim that 
there is a vertex $y_0 \in N^R_{\leq r-\alpha-2+t}(z)$  that lies at red-distance at least $2\alpha-2t+3$ from every vertex of $I\cup S$. To see this, first note that any shortest red path $P$ from $z$ to a vertex of $N^R_{r-\alpha-2+t}(z)$ contains at most $2(2\alpha-2t+2)+1$ vertices of $N^R_{\le 2\alpha-2t+2}(u)$ for every $u\in I\cup S$. This is because any set of more than $2(2\alpha-2t+2)+1$ vertices of $P$ contains a pair at red-distance more than $2(2\alpha-2t+2)$, whereas $N^R_{\le 2\alpha-2t+2}(u)$ has red diameter at most $2(2\alpha-2t+2)$.

Now since
\begin{align*}
|I\cup S|(2(2\alpha-2t+2)+1)
&\le
(\alpha+t)(2(2\alpha-2t+2)+1)\\
&=
4(\alpha^2-t^2)+5\alpha+5t\\
&=r-\alpha-2+t-4t(t-1)\\
&\leq
r-\alpha-2+t=|V(P)|-1,
\end{align*}
there exists a vertex $y_0$ as claimed.

Note that $S\cup\{y_0\}$ fulfills~\ref{g1} and~\ref{g2}, and thus, by our choice of~$t$, either~\ref{g3} fails for $S\cup\{y_0\}$, or $t=\alpha$ (and thus $t=|I|$ and $\{\ell(s):s\in S\}=I$). In either case, 
$\ell(y_0)=\ell(s)$ for some $s\in S$.  Note that $s$ is blue-adjacent to $y_0$ because of \eqref{eq:samelabels} and the fact that $y_0$ has no red neighbors in $S$.  Also note that $y_0$ has no other blue neighbors in $S$ because of \ref{g3} and \eqref{eq:differentlabels}.  Hence $N^B(y_0)\cap S=\{s\}$.

Let $T\subseteq S$ be the set of endpoints of all
alternating paths of the form $y_0s_1\ldots y_{k-1}s_k$, where for each $1\le i\leq k$,
$s_i\in S$, $s_iy_i$ is red and 
$y_{i-1}s_i$ is blue, and hence $\ell(y_{i-1})=\ell(s_i)$ by~\eqref{eq:differentlabels}. Note $s_1=s$, and by~\ref{g1} and our choice of $y_0$, the vertices of 
all such paths are contained in $N^R_{\leq r-\alpha-1+t}(z)$.

We claim that the set $\ell(N^R(T))=N^B(N^R(T))\cap I$ of labels of the vertices in
$N^R(T)$ satisfies 
\begin{equation}\label{eq:labels}
\ell(N^R(T))\subseteq\ell(T). 
\end{equation}
Indeed, suppose on
the contrary that we have an alternating path $y_0s_1\ldots y_{k-1}s_k$ with $s_k\in T$ of the form described above and suppose that $s_k$ has a red neighbour $y$ with $\ell(y)\notin \ell(T)$. If
there exists $s\in S\setminus T$ with $\ell(s)=\ell(y)$ then by~\eqref{eq:samelabels}, $ys$ is an edge,
which by \ref{g2} is blue, and so
the alternating path $y_0s_1\ldots y_{k-1}s_kys$ shows that $s$ should have
been in $T$. If there is no $s\in S$ with label $\ell(y)$ (i.e. $\ell(y)\not\in \ell(S)$), then
$t<|I|\leq\alpha$, and the set
$S':=(S\setminus\{s_1,\ldots,s_k\})\cup\{y_0,y_1,\ldots,y_{k-1},y\}$ has the
following properties:
\begin{enumerate}
\item $|S'|=t+1,$
\item $S'\subseteq N^R_{\leq r-\alpha-1+t}(z)$,
  \item every $s\in S'$ is at red-distance at least $2\alpha-2t+1=2\alpha-2(t+1)+3$ from
    every vertex of $I\cup (S'\setminus\{s\})$, 
  \item $\ell(\{y_0,\ldots,y_{k-1}\})=\ell(\{s_1,\ldots,s_k\})$, and
    $\ell(y)\notin\ell(S)=\ell(S')\setminus\{y\}$. % by assumption.
\end{enumerate}
Therefore the set $S'$ is $(t+1)$-good, contradicting our choice
of $S$. This verifies \eqref{eq:labels}.

\begin{figure}[htbp]
\begin{center}
\includegraphics[width=5in]{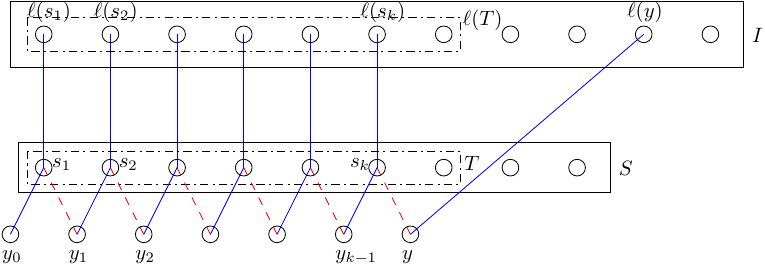}
\caption{Verifying \eqref{eq:labels}.}\label{fig:goodSet}
\end{center}
\end{figure}

To finish the proof, we define the graph $G'$ by removing the $|T|$
monochromatic blue subgraphs $\{N^B_{\leq 2}(v):v\in\ell(T)\}$ from $G$.  By \eqref{eq:labels}, the
independent set $T$ together with $N^R(T)$ is removed, and clearly
$N^B(T)$ is removed as well.  Hence $\alpha(G')\leq\alpha-|T|$. Then as
before, since by induction $G'$ has a covering set of at most
$\alpha(G')$ monochromatic subgraphs of diameter at most
$f(\alpha(G'))<f(\alpha)$, these together with $\{N^B_{\leq 2}(v):v\in\ell(T)\}$ cover $G$. Hence the proof is complete.

\section{Concluding Remarks}

The most immediate open problem arising from our work here is to improve our quadratic estimate $f(\alpha)\leq 8\alpha^2+12\alpha+6$ for the function $f$ introduced in Conjecture~\ref{con:DKMS}. As mentioned in the introduction, we are not aware of any nontrivial lower bound on $f(\alpha)$, and in particular whether or not it can be bounded above by an absolute constant independent of $\alpha$.

Another natural open problem is to consider edge colourings with more than two colours. This question was also raised in~\cite{DKMS}, based on the following old and notoriously difficult open problem known as Ryser's Conjecture.

\begin{conjecture}[Ryser]\label{con:Ryser}
For every $r$-partite hypergraph, $\tau(H)\leq (r-1)\nu(H)$.

Equivalently, for every graph $G$ and every $r$-colouring of the edges of $G$, there exists a covering set for $G$ of at most $(r-1)\alpha(G)$ monochromatic components.
\end{conjecture}

When $r=2$, this is K\H onig's theorem (Theorem \ref{thm:K1}/\ref{thm:K2}), and the $r=3$ case was proved by Aharoni \cite{A}. While Ryser's conjecture is still open for all $r\geq 4$, Tuza \cite{T1, T2} proved that the conjecture holds for $4\leq r\leq 5$ when $\nu(H)=1$ (equivalently $\alpha(G)=1$; i.e.~$G$ is complete).

Mili\'cevi\'c \cite{M1, M2} proved  a bounded diameter strengthing of Ryser's conjecture when $\alpha(G)=1$ and $3\leq r\leq 4$.  He also conjectured a bounded diameter strengthening of Ryser's conjecture in the case $\alpha(G)=1$.  The bounds of Mili\'cevi\'c were improved in \cite{DKMS}  and Mili\'cevi\'c's conjecture was generalized further.

\begin{conjecture}[DeBiasio, Kamel, McCourt, Sheats \cite{DKMS}]\label{con:DKMS2}
There exists a function $f:\NN\times \NN\to \NN$ such that for every graph $G$ and every $r$-colouring of the edges of $G$, there exists a covering set for $G$ of at most $(r-1)\alpha(G)$ monochromatic subgraphs, each of diameter at most $f(r,\alpha(G))$.
\end{conjecture}

Despite the fact that Ryser's conjecture is still open, it is conceivable that one could prove Conjecture \ref{con:DKMS2} under the assumption that Ryser's conjecture holds.

\paragraph{Acknowledgement.} This project was initiated at the Second Graph Theory in the Andes Workshop in March 2024. The authors wish to thank the Center for Mathematical Modeling  for its support sponsoring this event  through ANID Basal Grant CMM FB210005.

We thank the referees for their careful reading of the paper.

\bibliographystyle{abbrv}
\bibliography{references}

\end{document}